\documentclass[11pt,twoside]{article}
\usepackage[latin9]{inputenc}
\usepackage[a4paper]{geometry}
\usepackage{array}
\usepackage{float}
\usepackage{multirow}
\usepackage{amsmath}
\usepackage{amsthm}
\usepackage{amssymb}
\usepackage{graphicx}
\usepackage{esint}

\makeatletter

\providecommand{\tabularnewline}{\\}

\@ifundefined{showcaptionsetup}{}{%
 \PassOptionsToPackage{caption=false}{subfig}}
\usepackage{subfig}
\makeatother

\begin{document}
\title{Performance of nonconforming spectral element method for Stokes problems}
\author{N. Kishore Kumar\thanks{BITS-Pilani Hyderabad Campus, Hyderabad, Email: naraparaju@hyderabad.bits-pilani.ac.in}$^{\star}$
and Shubhashree Mohapatra\thanks{IIIT Delhi, Email: subhashree@iiitd.ac.in}$^{\dagger}$\\
 }
\date{~}
\maketitle
\begin{abstract}
In this paper, we study the performance of the non-conforming least-squares
spectral element method for Stokes problem. Generalized Stokes problem
has been considered and the method is shown to be exponential accurate.
The numerical method is nonconforming and higher order spectral element
functions are used. The same order spectral element functions are
used for both velocity and pressure variables. The normal equations
in the least-squares formulation are solved efficiently using preconditioned
conjugate gradient method. Various test cases are considered including
the Stokes problem on curvilinear domains, Stokes problem with mixed
boundary conditions and a generalized stokes problem in $\mathbb{R}^{3}$
to verify the accuracy of the method. 
\end{abstract}

\section{Introduction}

The stationary Stokes equations are linearization of the stationary
Navier-Stokes equations. These equations describe the flow of the
incompressible fluids. Stokes problems arise in various applications
in physics and engineering. The numerical solution of stationary Stokes
equations has been extensively studied in the literature. Finite difference
method, Finite element method, Spectral methods, discontinuous Galerkin
methods, Least-squares finite element methods and meshless methods
etc. are popular methods among the numerical techniques to solve the
Stokes problem.

Standard central differences do not give stable discretizations on
uniform grids to Stokes equations due to the failure of discrete inf-sup
condition \cite{CHEN}. The marker and cell method (MAC Scheme) is
simple and more efficient numerical scheme for solving Stokes and
Navier-Stokes problems. Finite difference MAC scheme has been studied
in \cite{CHEN,itoqiao}. Details of the other finite difference schemes
can be found in \cite{songli,srtrickwerda1,srtrickwerda2}. Finite
volume MAC scheme has been studied in \cite{ruili} and the references
therein.

Finite element method for Stokes problem has been widely studied.
Standard Galerkin formulation is viewed as saddle-point problem. Mixed
FEM impose restrictions such as inf-sup or Ladyzhenskaya-Babuska-Brezzi
(LBB) condition while choosing approximation spaces for different
unknowns, because of which same order polynomials can not be used
for different unknowns. So it uses two different finite element spaces
for velocity and pressure respectively. The analysis of mixed FEM
is based on the theory of saddle point problem which has been developed
in \cite{babuska,brezzi}. Use of different order polynomial spaces
makes computation cumbersome. Also choosing a suitable pair of polynomial
order spaces is not easy in general. Mass conservation is an another
issue in the approximation of the solution of the Stokes problem. 

To overcome the problem of using different finite element spaces for
velocity and pressure variables, in stabilized finite element method,
the standard bilinear form is modified such that any pair of finite
element spaces can be chosen. Stabilized finite element method for
Stokes problem has been introduced in \cite{brezzidoug,hugesfranka2,hugesfranka1}.
Details of several other stabilized finite element formulations can
be found in \cite{arnoldbrezzi,bar,BLA1,bofbrezzi,cod,douglas} and
the references therein. Divergence-free finite element methods have
been studied in \cite{blank,crouz,muye,wangwang} and the references
therein. Divergence-free methods eliminates pressure variable from
the saddle point system and results in positive definite linear systems
and also avoids the mass conservation issue.

Least squares methods have gained lots of attention for solving differential
equations over the last few decades \cite{AKS,BD,EASON,gunzbochev,JN1}.
These methods offer various advantages compared to Galerkin methods
(known as mixed methods), when applied to system of differential equations.
First of all, least squares based methods are free from any stabilizing
condition or parameters as mentioned above. Apart from this, least
squares methods always lead to a symmetric positive definite system
of equations when applied to linear problems which is one of the most
desirable property while solving a system of equations. In the least-squares
formulation Stokes problem is transformed into a first-order system.
There are different first order formulations for Stokes system. For
example velocity-vorticity-pressure formulation \cite{AMARACHA,BochevGunzburger,bramblepasc,CA1,CY,DU1,dubois,JN1},
velocity-stress-pressure formulation \cite{BD1,bramblepasc,kimshin}
and acceleration-pressure formulation \cite{chang}. 

Though least squares methods offer various theoretical and computational
advantages, it makes users skeptical when mass conservation is measured.
Continuity equation is minimized in a least squares sense which raises
this issue \cite{DG}. Here we briefly discuss existing literature
addressing mass conserving properties of least squares methods. Chang
et. al. \cite{chang2} have proposed a least squares finite element
method, which enforces continuity using Lagrange multiplier technique
in the least squares functional. The proposed method works well while
conserving mass properly, however looses one of the desirable property,
the positive definiteness of the system. A comparative study for three
different formulations of Stokes equations has been provided by Bolton
et. al. in \cite{BOTH}. The investigation concludes that lack of
mass conservation may lead to extremely poor results. Least squares
method for velocity-vorticity-pressure formulation has been used by
Proot et. al. \cite{PROO3}. An interesting observation from this
investigation is: least squares methods might lack good mass conservation
property, but provide better momentum conservation. Applying weights
to continuity equation in the least squares functional might disturb
the momentum conservation.

Stokes problem also has been studied using spectral/spectral element
methods. Spectral method has been proposed for this problem in \cite{schumack}.
Proot et. al. \cite{PROO1,PROO2,PROO3} have proposed a least squares
spectral element scheme for Stokes equations in velocity-vorticity-pressure
formulation. Numerical results conclude that the pressure variable
to be one order less accurate compared to velocity variable as pressure
is not prescribed on the boundary. Least-squares spectral collocation
method have been studied in \cite{HEIN,K1M1}. If the data in the
given problem is analytic then these methods gives exponential convergence.

Nonconforming methods like mortar finite element methods \cite{belgacem,belgacem1},
discontinuous Galerkin methods \cite{bursta,cockkan,lisunyang,montfer,monfer1}
and some other in \cite{apelkempf,BOCHEVLAI,BLO1,BLO2} provide numerical
approximation to the Stokes equations. Meshless methods always avoid
the problem of mesh generation on the domains with complex geometries.
Meshless methods for the Stokes problem have been studied in \cite{ahmsiraj,desimurq,xli,tanzhang,traskmaxi}.

In this article we have also considered the generalized Stokes problem.
This problem occurs in the numerical treatment of the time dependent
Navier-Stokes equations. Generalized Stokes equations looks very similar
to the Stokes equations. This problem has been studied in \cite{bar,bur,but,callam,chou,cod,lar1,nafa,sarin}.

In \cite{mahapravir} an exponentially accurate non-conforming least-squares
spectral element method for Stokes equations on non-smooth domains
has been proposed. This is different from the standard least-squares
FEM formulations where the Stokes system is converted into a first
order system as mentioned earlier. The minimizing functional in the
least-squares formulation includes the residuals in the partial differential
equations and residuals in the boundary conditions in appropriate
Sobolev norms. The method is nonconforming and higher order spectral
element functions have been used. 

In this article we study the performance of this method for Stokes
problem on smooth domains including curvilinear domains with different
boundary conditions and also for the generalized Stokes problem. The
normal equations in the least-squares formulation are solved using
preconditioned conjugate gradient method without storing the matrix.
Exponential accuracy of the method is verified through various numerical
tests.

This paper is organized as follows. Section 2, introduces few notations
that are needed for our analysis and proposes stability estimates
for velocity and pressure components in generalized Stokes equations.
The numerical scheme and error estimates are described in Section
3. Numerical results are presented in Section 4. Finally we conclude
with Section 5.

Here we give some notations and define required function spaces. Let
$\Omega\in\mathbb{R}^{2},$ be an open bounded set with sufficiently
smooth boundary $\partial\Omega$. $H^{m}(\Omega)$ denotes the Sobolev
space of functions with square integrable derivatives of integer order
less than or equal to $m$ on $\Omega$ equipped with the norm 
\[
\left\Vert u\right\Vert _{H^{m}(\Omega)}^{2}=\sum_{\left|\alpha\right|\leq m}\left\Vert D^{\alpha}u\right\Vert _{L^{2}(\Omega)}^{2}.
\]
Further, let $I=(-1,1).$ Then we define fractional norms $(0<s<1)$
by : 
\begin{eqnarray*}
\left\Vert w\right\Vert _{s,I}^{2}=\|w\|_{0,I}^{2}+\int_{I}\int_{I}\frac{\left|w(\xi)-w(\xi^{\prime})\right|^{2}}{\left|\xi-\xi^{\prime}\right|^{1+2s}}d\xi d\xi^{\prime},
\end{eqnarray*}
where $I$ denotes an interval contained in $\mathbb{R}.$ Moreover,
\begin{align*}
\|w\|_{1+s,I}^{2}=\|w\|_{0,I}^{2}+\left\Vert \:\frac{\partial w}{\partial\xi}\:\right\Vert _{s,I}^{2}+\left\Vert \:\frac{\partial w}{\partial\eta}\:\right\Vert _{s,I}^{2}\;.
\end{align*}
We shall denote the vectors by bold letters. For example, $\mathbf{u}=(u_{1},u_{2})^{T}$,
$\mathbf{H}^{k}(\Omega)=H^{k}(\Omega)\times H^{k}(\Omega),$ etc.
The norms are given by $||\mathbf{u}||_{k,\Omega}^{2}=||u_{1}||_{k,\Omega}^{2}+||u_{2}||_{k,\Omega}^{2}$
for $\mathbf{u}\in\mathbf{H}^{k}(\Omega),$ $\left\Vert \mathbf{u}\right\Vert _{s,I}^{2}=||u_{1}||_{s,I}^{2}+||u_{2}||_{s,I}^{2},$
etc. 

\section{Discretization and Stability Estimate}

In this section we describe the discretization of the domain and derive
the numerical formulation. 

\subsection{Generalized Stokes problem }

Consider the generalized Stokes equations in $\Omega\subset\mathbb{R}^{2}$,
with sufficiently smooth boundary $\partial\Omega=\Gamma$ (as shown
in fig. 1). 
\begin{align}
\alpha\mathbf{u}-\nu\Delta\mathbf{u}+\nabla p & =\mathbf{f}\quad\mbox{in}\quad\Omega\label{eqn2.21}\\
-\nabla\cdot\mathbf{u} & =h\quad\mbox{in}\quad\Omega\label{eqn2.22}\\
\mathbf{u} & =\mathbf{g}\quad\mbox{on}\quad\partial\Omega.\label{eqn2.23}
\end{align}
Here, $\mathbf{u}$ is the velocity field, $p$ is the pressure and
$\alpha,\nu>0$. Assume that the positive parameters $\alpha,\nu$
are not simultaneously zero. When $\alpha=0$ it reduces to Stokes
problem.

Let $h\in L^{2}(\Omega)$ and such that $\int_{\Omega}\,h\,dx=0$
and $\mathbf{g}$ satisfies the compatibility condition $\int_{\Gamma}\,\mathbf{g}.\mathbf{n}=0$
where $\mathbf{n}$ is the unit outward normal to $\Gamma.$ Let $\mathbf{f}\in\mathbf{L}^{2}(\Omega),\mathbf{g}\in\mathbf{H}^{\frac{1}{2}}(\partial\Omega).$
Then the generalized Stokes problem (1-3) has the solution $(\mathbf{u},p)\in\mathbf{H}^{1}(\Omega)\times L^{2}(\Omega),$
where $p$ is unique up to a additive constant. $p$ can be obtained
uniquely in $L^{2}(\Omega)/\mathbb{R}$ or in $L_{0}^{2}(\Omega)=\{u\in L^{2}(\Omega)|\int_{\Omega}u=0\}.$ 

Further, the following regularity estimate also holds good. This fundamental
regularity estimate is based on ADN (Agmin-Douglis-Nirenberg) theory
\cite{ADN}.\medskip{}

\subsubsection*{Regularity Estimate}

Let $\Omega$ be an open bounded subset of class $C^{r},r=\textrm{max}(m+2,2)$.
For $\mathbf{u}\in\mathbf{W}^{1,2}(\Omega)$, $p\in L^{2}(\Omega)$
being solutions of the generalized Stokes equations (\ref{eqn2.21})-(\ref{eqn2.23})
and for $\mathbf{f}\in\mathbf{W}^{m,2}(\Omega)$, $h\in W^{m+1,2}(\Omega)$
and $\mathbf{g}\in\mathbf{W}^{m+\frac{3}{2},2}(\Gamma)$, then $\mathbf{u}\in\mathbf{W}^{m+2,2}(\Omega)$,
$p\in W^{m+1,2}(\Omega)$ and there exists a constant $C_{0}(\alpha,m,\Omega)$
such that 
\begin{eqnarray}
||u||_{\mathbf{W}^{m+2,2}(\Omega)}+||p||_{W^{m+1,2}(\Omega)/\mathbb{R}}\leq C_{0}\left(||\mathbf{f}||_{\mathbf{W}^{m,2}(\Omega)}+||h||_{W^{m+1,2}(\Omega)}+||\mathbf{g}||_{\mathbf{W}^{m+\frac{3}{2},2}(\Gamma)}\right) &  & \blacksquare\label{eq:}
\end{eqnarray}
Let $\mathcal{L}(\mathbf{u},p)$ and $\mathcal{D}\mathbf{u}$ be the
differential operators for the momentum equations and the continuity
equation respectively. Thus, 
\begin{align*}
\mathcal{L}(\mathbf{u},p)=\alpha\mathbf{u}-\nu\Delta\mathbf{u}+\nabla p
\end{align*}
and 
\begin{align*}
\mathcal{D}(\mathbf{u})=-\nabla\cdot\mathbf{u}.
\end{align*}

\subsection{Discretization and spectral element functions}

\begin{figure}[H]

~~~~~~~~~~\includegraphics[scale=0.8]{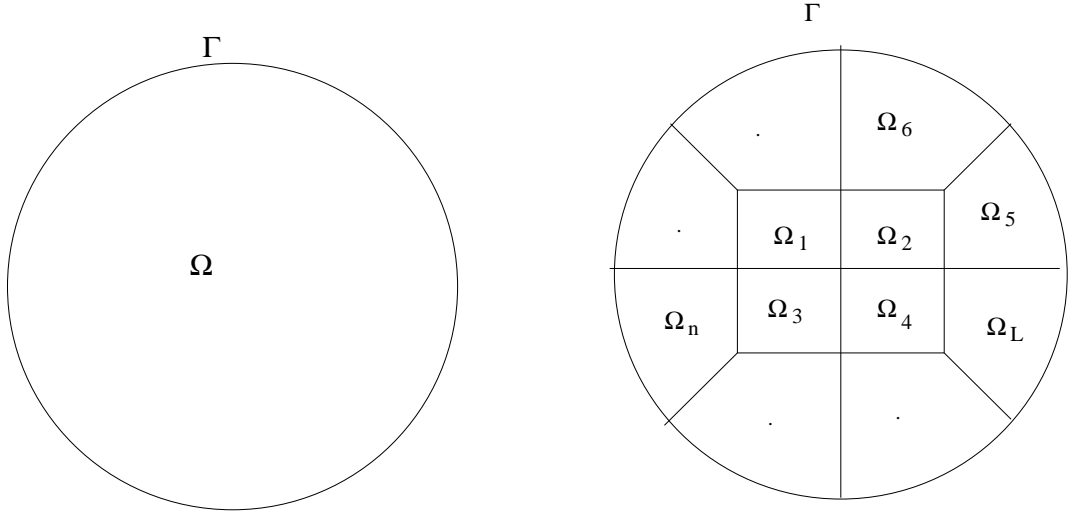}

\caption{Domain $\Omega$ and its discretization}
\end{figure}

The domain $\Omega$ is divided into $L$ quadrilaterals $\Omega_{1},\Omega_{2}...\Omega_{L}$
as shown in the figure 1 (some of them are curvilinear elements).
Let $\boldsymbol{x}=(x_{1},x_{2})$ be any point in the domain. A
set of nonconforming spectral element functions are defined on these
elements which are a sum of tensor products of polynomials of degree
$W$. Let $S$ denote the master element $S=(-1,1)^{2}$. Now there
is an analytic map $M_{l}(\xi,\eta)$ from $S$ to $\Omega_{l}$ which
has an analytic inverse (blending function mapping \cite{gordanhall})
\begin{eqnarray*}
x_{1}=X_{1}^{l}(\xi,\eta)\,\,\textrm{and}\,\,x_{2}=X_{2}^{l}(\xi,\eta).
\end{eqnarray*}

Define the spectral element functions $\hat{\mathbf{u}}_{l}$ and
$\hat{p}_{l}$ on $S$ by 
\begin{align*}
\hat{\mathbf{u}}(\xi,\eta)=\sum_{i=0}^{W}\sum_{j=0}^{W}\mathbf{a}_{i,j}\:\xi^{i}\eta^{j},\quad\hat{p}(\xi,\eta)=\sum_{i=0}^{W}\sum_{j=0}^{W}b_{i,j}\xi^{i}\eta^{j}.
\end{align*}

Then $\mathbf{u}_{l}$ and $p_{l}$ on $\Omega_{l}$ are given by
\begin{eqnarray*}
\mathbf{u}_{l}(x_{1},x_{2})=\hat{\mathbf{u}}_{l}(M_{l}^{-1})\,\,\textrm{and}\,\,p_{l}(x_{1},x_{2})=\hat{p}(M_{l}^{-1}).
\end{eqnarray*}
Let $\Pi^{L,W}=\left\{ \left\{ \mathbf{u}_{l}\right\} _{1\leq l\leq L},\left\{ p_{l}\right\} _{1\leq l\leq L}\right\} $
be the space of spectral element functions consisting of the above
tensor products of polynomials of degree $W$.

\subsection{Stability estimate}

Let $J_{l}(\xi,\eta)$ be the Jacobian of the mapping $M_{l}(\xi,\eta)$
from $S=(-1,1)^{2}$ to $\Omega_{l}$ for $l=1,2,...,L.$ Now 
\[
\intop_{\Omega_{l}}\left|\mathcal{L}\mathbf{u}_{l}\right|^{2}dx_{1}dx_{2}=\intop_{S}\left|\mathcal{L}\hat{\mathbf{u}}_{l}\right|^{2}J_{l}\,d\xi d\eta.
\]
Define $\mathcal{L}_{l}\hat{\mathbf{u}}_{l}=\mathcal{L}\hat{\mathbf{u}}_{l}\sqrt{J_{l}}.$
Then 
\[
\intop_{\Omega_{l}}\left|\mathcal{L}\mathbf{u}_{l}\right|^{2}dx_{1}dx_{2}=\intop_{S}\left|\mathcal{L}_{l}\hat{\mathbf{u}}_{l}\right|^{2}\,d\xi d\eta.
\]
Similarly, we define $\mathcal{D}_{l}\hat{\mathbf{u}}_{l}=\mathcal{D}\hat{\mathbf{u}}_{l}\,\sqrt{J_{l}}.$ 

Since the approximation is nonconforming, to enforce the continuity
along the inter element boundaries we introduce the jumps in $\mathbf{u},\mathbf{u}_{x_{1}},\mathbf{u}_{x_{2}}$
and $p$ in suitable Sobolev norms. Let the edge $\gamma_{s}$ be
common to the adjacent elements $\Omega_{l}\,\,\textrm{and}\,\,\Omega_{m}.$
Assume that edge $\gamma_{s}$ is the image of $\eta=1$ under the
map $M_{l}$ which maps $S$ to $\Omega_{l}$ and also the image of
$\eta=-1$ under the map $M_{m}$ which maps $S$ to $\Omega_{m}$.
Then the jumps along the inter-element boundaries are defined as
\begin{align*}
 & \left\Vert [\mathbf{u}]\right\Vert _{0,\gamma_{s}}^{2}=\left\Vert \hat{\mathbf{u}}_{m}(\xi,-1)-\hat{\mathbf{u}}_{l}(\xi,1)\right\Vert _{0,I}^{2},\\
 & \left\Vert [\mathbf{u}_{x_{k}}]\right\Vert _{\frac{1}{2},\gamma_{s}}^{2}=\|({\hat{\mathbf{u}}_{m})}_{x_{k}}(\xi,-1)-{(\hat{\mathbf{u}}_{l})}_{x_{k}}(\xi,1)\|_{\frac{1}{2},I}^{2},\\
 & \left\Vert [p]\right\Vert _{\frac{1}{2},\gamma_{s}}^{2}=\left\Vert \hat{p}_{m}(\xi,-1)-\hat{p}_{l}(\xi,1)\right\Vert _{\frac{1}{2},I}^{2}.
\end{align*}
Here $I=(-1,1).$ The expressions on the right hand side in the above
equation are given in the transformed coordinates $\xi$ and $\eta.$ 

Let us consider the boundary condition. Let $\gamma_{s}\subseteq\partial\Omega\cap\Omega_{l}$
be the image of $\xi=1$ under the mapping $M_{l}$ which maps $S$
to $\Omega_{l}.$ Then
\begin{eqnarray*}
||\mathbf{u}_{l}||_{\frac{3}{2},\gamma_{s}}^{2}=||\hat{\mathbf{u}}_{l}(1,\eta)||_{\frac{3}{2},I}^{2}.
\end{eqnarray*}

Let $\mathbf{u},p\in\Pi^{L,W}$. We now define the quadratic form
\begin{align}
\mathcal{V}^{L,W}(\mathbf{u},p) & =\sum_{l=1}^{L}||\mathcal{L}(\mathbf{u}_{l},p_{l})||_{0,\Omega_{l}}^{2}+\sum_{l=1}^{L}||\mathcal{D}\mathbf{u}_{l}||_{1,\Omega_{l}}^{2}\nonumber \\
 & +\sum_{\gamma_{s}\subseteq\bar{\Omega}\setminus\partial\Omega}\left(||[\mathbf{u}]||_{0,\gamma_{s}}^{2}+\sum_{k=1}^{2}||[\mathbf{u}_{x_{k}}]||_{\frac{1}{2},\gamma_{s}}^{2}+||[p]||_{\frac{1}{2},\gamma_{s}}^{2}\right)+\sum_{\gamma_{s}\subseteq\partial\Omega\cap\Omega_{l}}||\mathbf{u}_{l}||_{\frac{3}{2},\gamma_{s}}^{2}.
\end{align}
Next, we define the quadratic form 
\begin{align}
\mathcal{U}^{L,W}(\mathbf{u},p)=\sum_{l=1}^{L}||\mathbf{u}_{l}||_{2,S}^{2}+\sum_{l=1}^{L}||p_{l}||_{1,S}^{2}.\label{eq3.4d}
\end{align}
Then, we have the following result. \\
\\
\textbf{Theorem 2.1:} For $W$ large enough there exists a constant
$C>0$ such that the estimate 
\begin{align}
\mathcal{U}^{L,W}(\mathbf{u},p)\leq C(\ln W)^{2}\mathcal{V}^{L,W}(\mathbf{u},p)\label{eqn101}
\end{align}
holds.

The proof of this one is very similar to Theorem 4.1 in \cite{mahapravir}. 

\section{Numerical scheme and error estimates}

In this section we describe the numerical scheme which is based on
the stability estimate Theorem 2.1 and state the error estimate. We
also brief the residual computations. 

\subsection{Numerical scheme}

We search for a minimizer which minimizes the sum of the residuals
in the partial differential equations, the residuals in the boundary
conditions and the jumps in the velocity variable, derivatives of
velocity variable and pressure variable along the inter element boundaries
in appropriate fractional Sobolev norms. To formulate the numerical
scheme we define a functional $\mathcal{R}^{L,W}(\mathbf{u},p)$,
closely related to the quadratic form $\mathcal{V}^{L,W}(\mathbf{u},p)$.

As defined in section 2, let $M_{l}(\xi,\eta)$ be a mapping from
$S$ to $\Omega_{l}.$ Let $\hat{\boldsymbol{x}}=(\xi,\eta).$ Let
$\mathbf{f}_{l}(\hat{\boldsymbol{x}})=\mathbf{f}(M_{l}(\xi,\eta))$,
$h_{l}(\hat{\boldsymbol{x}})=h(M_{l}(\xi,\eta))$, for $l=1,2,\ldots,L.$
Let $J_{l}(\hat{\boldsymbol{x}})$ denote the Jacobian of the mapping
$M_{l}$ from $S$ to $\Omega_{l}.$ Define $\mathbf{F}_{l}(\hat{\boldsymbol{x}})=\mathbf{f}_{l}(\hat{\boldsymbol{x}})\sqrt{J_{l}(\hat{\boldsymbol{x}})}$,
$H_{l}(\hat{\boldsymbol{x}})=h_{l}(\hat{\boldsymbol{x}})\sqrt{J_{l}(\hat{\boldsymbol{x}}).}$
Now consider the boundary condition $\mathbf{u}=\mathbf{g}$ on $\partial\Omega$.
Let $\gamma_{s}\subseteq\partial\Omega$ be the image of $\xi=1$
under the mapping $M_{l}$ which maps $S$ to $\Omega_{l}.$ Let $\mathbf{g}_{l}=\mathbf{g}(M_{l}(1,\eta))$.
We now define the least-squares functional 
\begin{align}
 & \mathcal{R}^{L,W}(\mathbf{u},p)=\sum_{l=1}^{L}\left\Vert \mathcal{L}_{l}\hat{\mathbf{u}}_{l}-\mathbf{F}_{l}\right\Vert _{0,S}^{2}+\sum_{l=1}^{L}||\mathcal{D}_{l}\mathbf{\hat{u}}_{l}-H_{l}||_{1,S}^{2}\nonumber \\
 & +\sum_{\gamma_{s}\subseteq\bar{\Omega}\setminus\partial\Omega}\left(||[\mathbf{u}]||_{0,\gamma_{s}}^{2}+\sum_{k=1}^{2}||[\mathbf{u}_{x_{k}}]||_{\frac{1}{2},\gamma_{s}}^{2}+||[p]||_{\frac{1}{2},\gamma_{s}}^{2}\right)+\sum_{\gamma_{s}\subseteq\partial\Omega\cap\Omega_{l}}||\mathbf{u}_{l}-\mathbf{g}_{l}||_{\frac{3}{2},\gamma_{s}}^{2}.
\end{align}
We now formulate the numerical scheme as follows: \textit{Find unique
($\mathbf{z},q)\in{\Pi}^{L,W}$ which minimizes the functional $\mathcal{R}^{L,W}(\mathbf{u},p)$
over all $(\mathbf{u},p)\in{\Pi}^{L,W}$. Here, ${\Pi}^{L,W}$ denotes
the space of spectral element functions.}

\subsubsection*{Error estimates}

\textbf{Theorem 3.1:} Let $(\mathbf{z},q)$ minimize $\mathcal{R}^{L,W}(\mathbf{u},p)$.
Then for $W$ large enough there exists constants $C$ and $b$ (being
independent of $W$) such that the estimate 
\begin{align}
\sum_{l=1}^{L}||\mathbf{z}_{l}(\hat{\boldsymbol{x}})-u_{l}(\hat{\boldsymbol{x}})||_{2,S}^{2}+\sum_{l=1}^{L}||q_{l}(\hat{\boldsymbol{x}})-p_{l}(\hat{\boldsymbol{x}})||_{1,S}^{2}\leq Ce^{-bW}
\end{align}
holds true.

Proof of the this theorem easily follows from theorem 4.2 from \cite{mahapravir}.
\\
\\
\textbf{\textit{Remark:}}\textit{ After obtaining a nonconforming
solution a set of corrections can be made such that velocity variable
$\mathbf{z}$ becomes conforming \cite{kishorekumar,CH}. So $\mathbf{z}\in\mathbf{H}^{1}(\Omega)$
and we have the following error estimate }
\begin{align}
||\mathbf{u}-\mathbf{z}||_{1,\Omega}+||p-q||_{0,\Omega}\leq Ce^{-bW}.
\end{align}

\subsection{Residue computations and preconditioner}

The solution is obtained at Gauss-Legendre-Lobatto (GLL) quadrature
points by minimizing the residue $\mathcal{R}^{L,W}(\mathbf{u},p).$
The normal equations obtained from the minimization are solved using
preconditioned conjugate gradient method (PCGM) without storing the
matrix. At each iteration step, PCGM requires only the action of the
matrix on a vector whose entries are the values of velocity and pressure
variable at GLL points and arranged in lexicographic order. The details
of the residual computations in each element and procedure of solving
the normal equations is shown in detail in \cite{DKU,kishorekumar,mahapravir}. 

Since the jumps across the inter-element boundaries in velocity and
pressure variables and also in the derivatives of velocity variable
are included in the numerical formulation, a small size of data have
to be interchanged between the elements at each iteration step of
the PCGM. The details of the preconditioner which we have used is
given below. 

Using Theorem 2.1, we have 
\begin{align}
\mathcal{U}^{L,W}(\mathbf{u},p)\leq C(\ln W)^{2}\mathcal{V}^{L,W}(\mathbf{u},p).\label{eqn10}
\end{align}
Using the trace theorem for Sobolev spaces, we get for some constant
$\tilde{C}$ 
\begin{align}
\tilde{C}\mathcal{V}^{L,W}(\mathbf{u},p)\leq\mathcal{U}^{L,W}(\mathbf{u},p)\leq C(\ln W)^{2}\mathcal{V}^{L,W}(\mathbf{u},p).
\end{align}
Hence the two quadratic forms $\mathcal{U}^{L,W}(\mathbf{u},p)$,
$\mathcal{V}^{L,W}(\mathbf{u},p)$ are spectrally equivalent and we
choose $\mathcal{U}^{L,W}(\mathbf{u},p)$ as our preconditioner so
that the condition number of the preconditioned system is polylogarithmic
in $W$, where $W$ is the degree of the polynomial. In each element
preconditioner consists of three blocks, first two blocks correspond
to $H^{2}$ norm of velocity variable and third block corresponds
to $H^{1}$ norm of pressure variable.

\section{Numerical results}

Here we verify the exponential convergence of the numerical scheme
by considering various numerical examples. The numerical examples
includes the Stokes equations on curvilinear domains, Stokes problem
with mixed boundary conditions and generalized Stokes equations in
two and three dimensions. Spectral element functions of higher order
of degree $W$ are used and uniform $W$ is used for all the elements
in the discretization. Let $\mathbf{z}$ and $q$ be the approximate
solutions of the velocity $\mathbf{u}$ and pressure $p$ respectively.
$\|E_{\mathbf{u}}\|_{1}=\frac{\left\Vert \mathbf{u}-\mathbf{z}\right\Vert _{1}}{\left\Vert \mathbf{u}\right\Vert _{1}}$
denotes the relative error in $\mathbf{u}$ in $H^{1}$ norm, $\|E_{p}\|_{0}=\frac{\left\Vert p-q\right\Vert _{0}}{\left\Vert p\right\Vert _{0}}$
denotes relative error in pressure in $L^{2}$ norm and ${\|E_{c}\|}_{0}$
denotes error in continuity equation in $L^{2}$ norm. 'Iter' denotes
the total number of iterations required to reach the desired accuracy.
In the case of Dirichlet boundary value problem, pressure is specified
to be zero at one point of the domain to ensure the uniqueness in
each problem.

\subsubsection*{Example 1: Generalized Stokes problem on $[0,1]^{2}$}

\begin{figure}[H]
~~~~~~~~~~~~~~~~~~~~~~~~~~~~~~~~~~~~~~\includegraphics[scale=0.7]{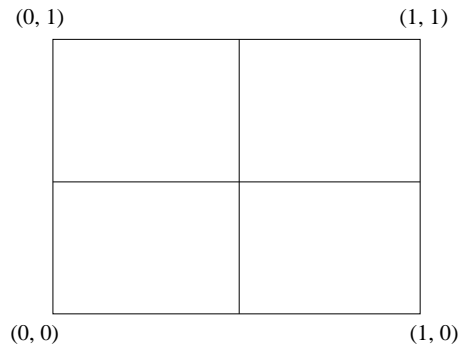}

\caption{Discretization of $[0,1]^{2}$}
\end{figure}

Consider the generalized Stokes equation (1-3) with $\alpha=1,\nu=1$
and $h=0$ on $[0,1]^{2}.$ Chosen the data such that 
\begin{eqnarray*}
u_{1}=sin\,\pi x_{1}\,\,sin\,\pi x_{2},\\
u_{2}=cos\,\pi x_{1}\,\,cos\,\pi x_{2},\\
p=150\,(x_{1}-\frac{1}{2})(x_{2}-\frac{1}{2})+c.
\end{eqnarray*}
The domain $[0,1]^{2}$ is divided into $4$ elements with equal step
size $h=\frac{1}{2}$ in both $x_{1}$ and $x_{2}$ directions (see
fig. 2). The approximate solution is obtained and the relative errors
$\|E_{\mathbf{u}}\|_{1}$ and $\|E_{p}\|_{0}$ for various values
of $W$ are shown in Table 1. Table 1 also shows the total number
of iterations of PCGM to reach the achieved accuracy and the error
in continuity equation $\|E_{c}\|_{0}.$ One can see that the errors
$\|E_{\mathbf{u}}\|_{1}$, $\|E_{p}\|_{0}$ and $\|E_{c}\|_{0}$ decays
exponentially. Fig. 3 shows the graph of log of the relative errors
$\|E_{\mathbf{u}}\|_{1}$ and $\|E_{p}\|_{0}$ against the log of
$W.$ The curves are almost linear. This shows the exponential decay
of the errors. 

\begin{table}[H]
~~~~~~~~~~~%
\begin{tabular}{ccccc}
\hline 
$W$  & $\|E_{\mathbf{u}}\|_{1}$  & $\|E_{p}\|_{0}$  & $\|E_{c}\|_{0}$ & Iter\tabularnewline
\hline 
2  & 5.139862909E-01  & 1.547329906E-01 & 1.177314788E-00 & 22\tabularnewline
3  & 9.320542876E-02  & 5.595518233E-02 & 1.701809578E-01 & 116\tabularnewline
4  & 1.644416626E-02  & 6.011788360E-03  & 3.465516273E-02 & 227\tabularnewline
5  & 1.611022198E-03 & 6.748736448E-04 & 3.739501208E-03 & 442\tabularnewline
6  & 1.693357558E-04 & 6.570180017E-05  & 4.172145413E-04 & 761\tabularnewline
7  & 2.471931274E-05  & 7.480811937E-06 & 6.553699066E-05 & 1204\tabularnewline
8  & 3.451649675E-06  & 4.728984810E-07  & 8.701982919E-06 & 1788\tabularnewline
\hline 
\end{tabular}

\caption{Error $\|E_{\mathbf{u}}\|_{1}$, $\|E_{p}\|_{0}$ and $\|E_{c}\|_{0}$
for different $W$}
\end{table}

\begin{figure}[H]
~~~~~~~~~\includegraphics[scale=0.85]{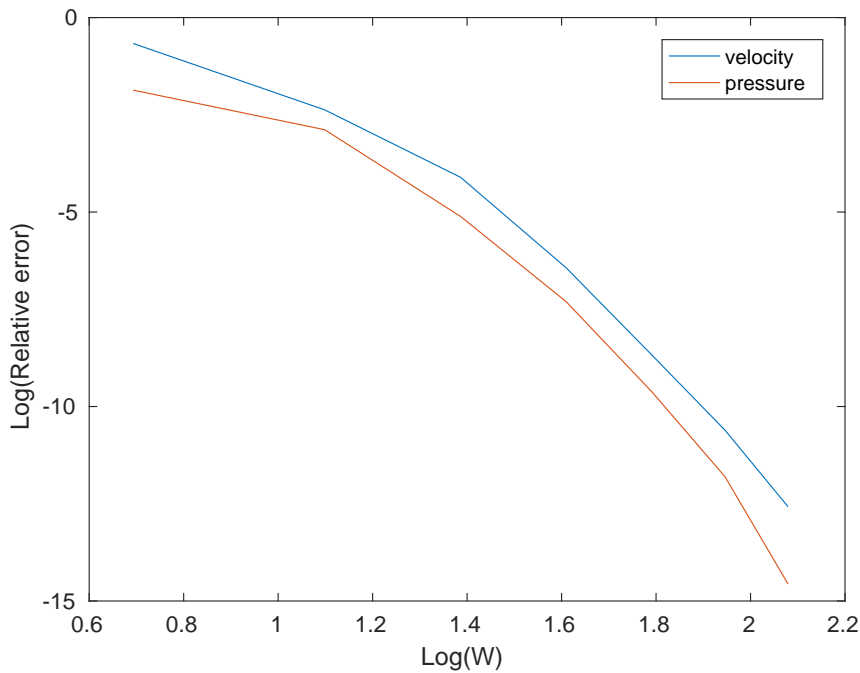}

\caption{Log of relative error vs. log($W$)}
\end{figure}

\subsubsection*{Example 2: Stokes problem involving Reynolds number}

Consider the following equation on $\Omega=[\frac{-1}{2},\frac{1}{2}]\times[0,1]$ 

\begin{eqnarray*}
-\frac{1}{Re}\triangle\mathbf{u}+\nabla p=\mathbf{f}\,\,\textrm{in}\,\,\Omega\\
\nabla.\mathbf{u}=0\,\,\textrm{in}\,\,\Omega\\
\mathbf{u=g}\,\,\textrm{on}\,\,\Gamma.
\end{eqnarray*}
 Choose the data such that 
\begin{eqnarray*}
u_{1}=1-e^{\lambda x_{1}}\,cos(2\pi x_{2}),\\
u_{2}=\frac{\lambda}{2\pi}\,e^{\lambda x_{1}}sin(2\pi x_{2}),\\
p=\frac{1}{2}e^{2\lambda x_{1}}+c.
\end{eqnarray*}
 Here $\lambda=\frac{Re}{2}-\sqrt{\frac{(Re)^{2}}{4}+4\pi^{2}}$ and
$Re$ is Reynolds number \cite{muye}. 

The domain is divided into 4 elements with step size $\frac{1}{2}$
in both directions. We have obtained the approximate solution of the
given Stokes system for different values of $Re=1,10,100,1000.$ Table
2 shows the relative errors $\|E_{\mathbf{u}}\|_{1}$ and $\|E_{p}\|_{0}$
for various values of $W$ for $Re=1,10.$ Table 3 shows the relative
errors $\|E_{\mathbf{u}}\|_{1}$ and $\|E_{p}\|_{0}$ for various
values of $W$ for $Re=100,1000.$ 

\begin{table}[H]
\begin{tabular}{cccc|ccc}
\hline 
 & $Re=1$ &  &  &  & $Re=10$ & \tabularnewline
\hline 
$W$ & $\|E_{\mathbf{u}}\|_{1}$ & $\|E_{p}\|_{0}$ & Iter & $\|E_{\mathbf{u}}\|_{1}$ & $\|E_{p}\|_{0}$ & Iter\tabularnewline
\hline 
4 & 2.366975497E-01 & 4.34268366E-01 & 12 & 1.131661632E-01 & 1.198261163E-01 & 40\tabularnewline
5 & 4.869971324E-02 & 1.76876076E-01 & 129 & \multirow{1}{*}{1.184006906E-02} & 1.227193197E-02 & 120\tabularnewline
6 & 8.623443013E-03 & 8.20262592E-02 & 321 & 3.004801200E-03 & 1.629441123E-03 & 195\tabularnewline
7 & 1.518313601E-03 & 1.30971235E-02 & 733 & 2.839200284E-04 & 3.512941780E-04 & 346\tabularnewline
8 & 2.667438050E-04 & 3.12678320E-03 & 1206 & 3.181384304E-05 & 4.256880275E-05 & 561\tabularnewline
9 & 4.391290170E-05 & 6.57675621E-04 & 1862 & 1.053751269E-06 & 9.059197820E-06 & 967\tabularnewline
10 & 7.408947632E-06 & 1.56727777E-04 & 2539 & 1.718964312E-07 & 2.683927048E-07 & 1843\tabularnewline
\hline 
\end{tabular}

\caption{Relative errors against $W$ for $Re=1,10$}

\end{table}

\begin{table}[H]
\begin{tabular}{cccc|ccc}
\hline 
 & $Re=100$ &  &  &  & $Re=1000$ & \tabularnewline
\hline 
$W$ & $\|E_{\mathbf{u}}\|_{1}$ & $\|E_{p}\|_{0}$ & Iter & $\|E_{\mathbf{u}}\|_{1}$ & $\|E_{p}\|_{0}$ & Iter\tabularnewline
\hline 
4 & 4.270224166E-01 & 1.971826154E-01 & 26 & 3.7944941768E-01 & 2.480889414E-01 & 36\tabularnewline
5 & 4.244262501E-02 & 1.618248455E-02 & 179 & 3.8779299147E-02 & 2.947725331E-02 & 654\tabularnewline
6 & 3.958119651E-03 & 1.747850562E-03 & 325 & 2.3764211892E-03 & 5.438172800E-04 & 1703\tabularnewline
7 & 3.201794986E-04 & 1.504988745E-04 & 514 & 1.9140466720E-04 & 5.837954856E-05 & 3190\tabularnewline
8 & 1.265052314E-05 & 3.527247772E-05 & 748 & 1.5032229108E-05 & 5.156661742E-06 & 4569\tabularnewline
9 & 3.187641118E-06 & 2.098243693E-06 & 936 & 3.5580252816E-06 & 1.602969000E-06 & 5721\tabularnewline
10 & 3.037312630E-07 & 1.211241028E-07 & 1271 & 7.451543690E-07 & 6.117219799E-07 & 7160\tabularnewline
\hline 
\end{tabular}

\caption{Relative errors against $W$ for $Re=100,1000$}
\end{table}

Fig. 4a shows the graph of log of the relative error $\|E_{\mathbf{u}}\|_{1}$
vs. log($W$) and fig. 4b shows the graph of log of the relative error
$\|E_{p}\|_{0}$ vs. log($W$) for $Re=1,10,100$ and $1000.$ Both
the graphs shows that the error decays exponentially. One can see
that the iteration count is high for $Re=1000$ to achieve the relative
errors of $O(10^{-7})$ and $O(10^{-6}),$ but the iteration count
is not high to achieve the accuracy of $O(10^{-4})$ or $O(10^{-5}).$ 

\begin{figure}[H]
\subfloat[Log($\|E_{\mathbf{u}}\|_{1}$) against log$(W)$]{\includegraphics[scale=0.6]{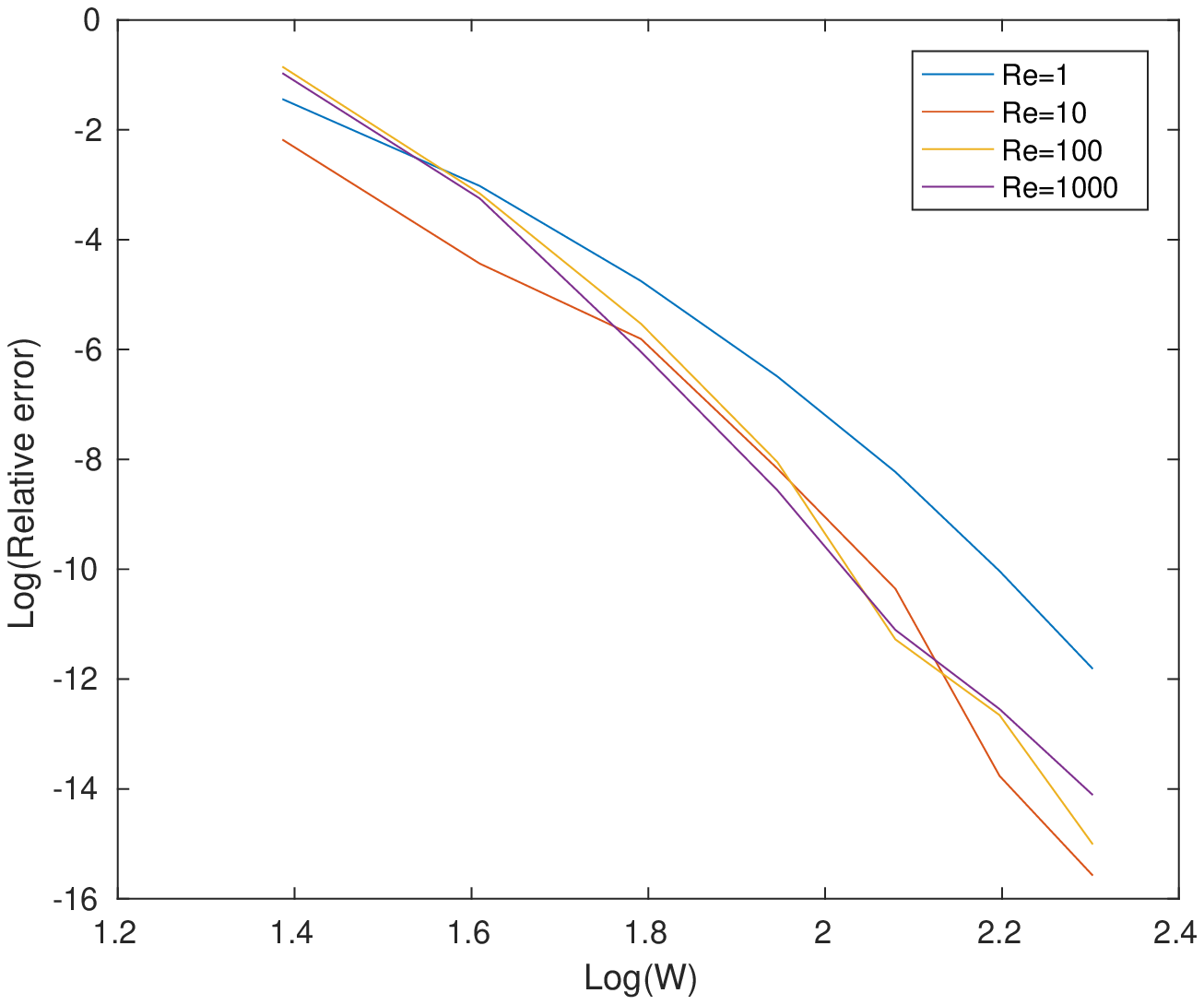}

}$\!\!\!\!\!\negthinspace\!$\subfloat[Log ($\|E_{p}\|_{0}$ ) against log($W$)]{\includegraphics[scale=0.64]{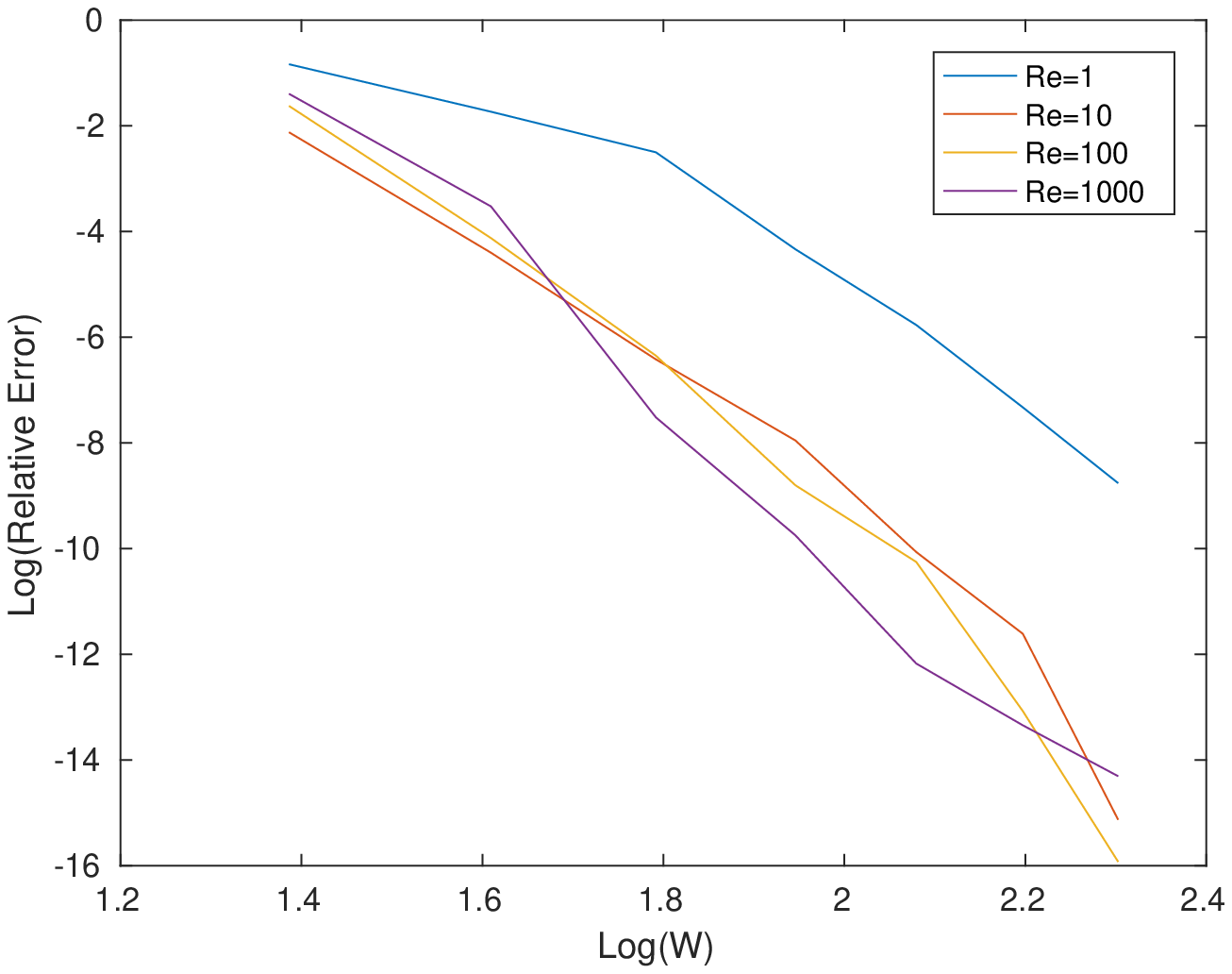}

}

\caption{Log of the relative errors against log($W$)}
\end{figure}

\subsubsection*{Example 3: Stokes problem on annular domain}

Consider the Stokes problem (problem (1-3) with $\alpha=0$, $\nu=1$
and $h=0$) on the annular domain $\Omega=\{(r,\theta):1\leq r\leq4\,\,\textrm{and}\,\,0\leq\theta\leq\frac{\pi}{2}\}$
with Dirichlet boundary condition on the boundary. 

\begin{figure}[H]
~~~~~~~~~~~~~~~~~~~~~~~~~~~~~~~~~~\includegraphics[scale=0.8]{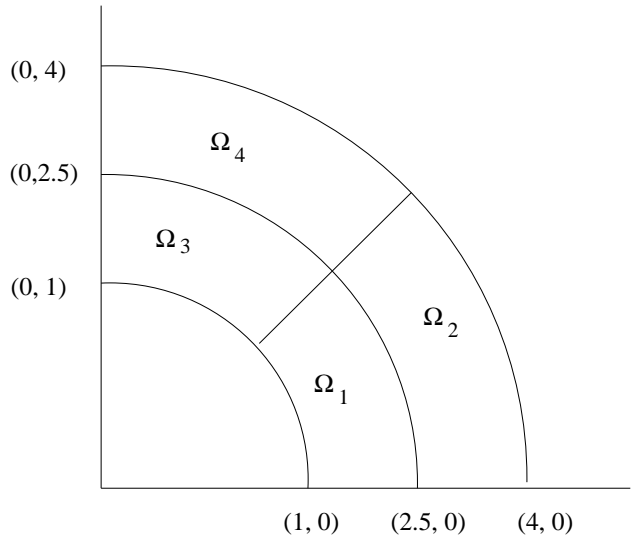}

\caption{Annular domain $\Omega$ and its discretization}

\end{figure}

The domain is divided into 4 curvilinear elements as shown in figure
5. Blending elements have been used \cite{gordanhall}. The data is
chosen such that 
\begin{eqnarray*}
u_{1}=20\,x_{1}\,x_{2}^{3},\\
u_{2}=5(x_{1}^{4}-x_{2}^{4}),\\
p=60x_{1}^{2}x_{2}-20x_{2}^{3}+c.
\end{eqnarray*}

Table 4 shows the relative errors $\|E_{\mathbf{u}}\|_{1}$ , $\|E_{p}\|_{0}$
and $\|E_{c}\|_{0}$ for various values of $W$. Fig. 6 shows the
log of the relative errors against log$(W).$ This shows that the
error decays exponentially in $\|E_{\mathbf{u}}\|_{1}$ and $\|E_{p}\|_{0}$
norms. 

\begin{table}[H]
~~~~~~~~~~~~~~~~~~~~~~~~~~~%
\begin{tabular}{ccccc}
\hline 
$W$ & $\|E_{\mathbf{u}}\|_{1}$ & $\|E_{p}\|_{0}$ & $\|E_{c}\|_{0}$ & Iter\tabularnewline
\hline 
2 & 3.73710182E-01 & 1.000929095E-00 & 558.6790137E-00 & 4\tabularnewline
3 & 5.47325472E-02 & 1.019245995E-01 & 149.7030088E-00 & 91\tabularnewline
4 & 9.27471997E-03 & 1.026087077E-02 & 17.92350320E-00 & 178\tabularnewline
5  & 2.27183319E-03 & 4.526209660E-03 & 7.016461256E-00 & 360\tabularnewline
6 & 4.61843438E-04 & 6.175899621E-04 & 1.262648344E-00 & 451\tabularnewline
7  & 4.52049220E-05 & 1.124195597E-04 & 1.507379559E-01 & 935\tabularnewline
8 & 8.47350232E-06 & 5.888507439E-06 & 2.122666602E-02 & 1252\tabularnewline
9 & 7.82707160E-07 & 1.396124323E-06 & 2.061182562E-03 & 1930\tabularnewline
10 & 1.52518872E-07 & 1.042122306E-07 & 3.131737628E-04 & 2748\tabularnewline
\hline 
\end{tabular}

\caption{$\|E_{\mathbf{u}}\|_{1}$, $\|E_{p}\|_{0}$ and $\|E_{c}\|_{0}$ for
various values of $W$}

\end{table}

\begin{figure}
~~~~~~~~~~~~\includegraphics[scale=0.85]{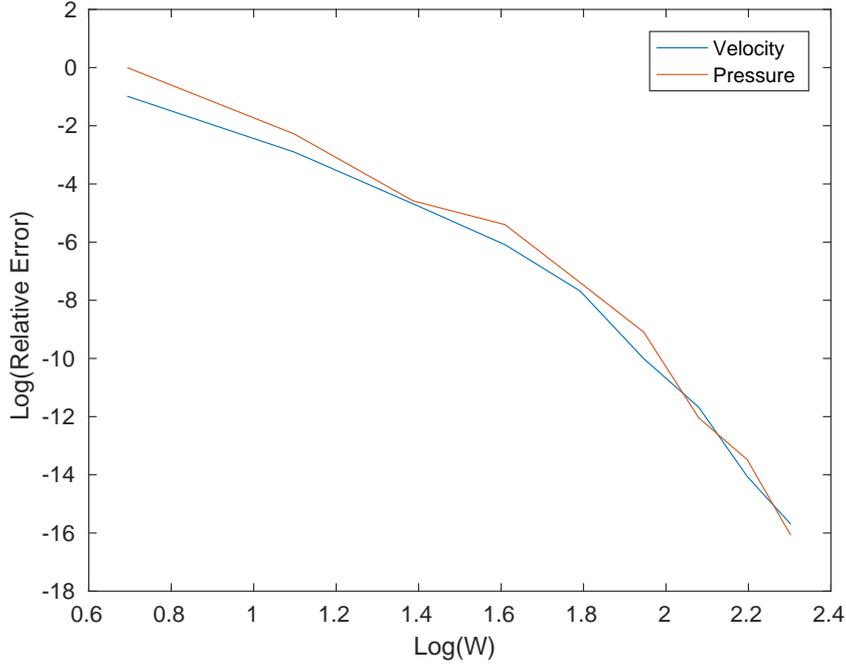}

\caption{Log of relative errors against log($W)$}

\end{figure}

\subsubsection*{Example 4. Stokes problem on a square domain with a circular hole}

Consider the Stokes problem on a square domain $[0,1]^{2}$ with a
circular hole where the circle is centered at $(0.5,0.5)$ with radius
$0.2$ (see fig. 7). As shown in the figure 7, we have the boundaries
on the sides of the unit square and also on the circle.

The data is chosen such that 
\begin{eqnarray*}
u_{1}=x_{1}+x_{2}^{2}-2x_{1}x_{2}+x_{1}^{3}-3x_{1}x_{2}^{2}+x_{2}x_{1}^{2},\\
u_{2}=-x_{2}-2x_{1}x_{2}+x_{2}^{2}-3x_{2}x_{1}^{2}+x_{1}^{3}-x_{1}x_{2}^{2},\\
p=x_{1}x_{2}+x_{1}+x_{2}+x_{1}^{3}x_{2}^{2}+c.
\end{eqnarray*}
\begin{figure}[H]

~~~~~~~~~~~~~~~~~~~~~~~~~~~~~\includegraphics[scale=0.8]{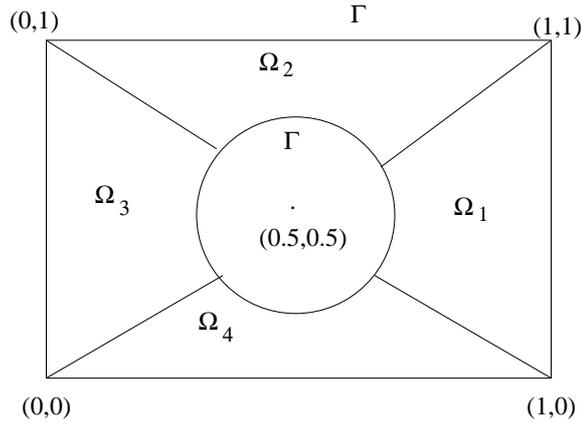}

\caption{The unit square with a circular hole}

\end{figure}

The given domain is decomposed into $4$ elements as shown in the
figure 7. Table 5 shows the relative errors $\|E_{\mathbf{u}}\|_{1}$
and $\|E_{p}\|_{0}$ for various values of $W$. 
\begin{table}[H]
~~~~~~~~~~~~~~~~~~~~~~%
\begin{tabular}{cccc}
\hline 
$W$ & $\|E_{\mathbf{u}}\|_{1}$ & $\|E_{p}\|_{0}$ & Iter\tabularnewline
\hline 
$2$ & 4.00582657E-01 & 9.966591533E-01 & 7\tabularnewline
\multicolumn{1}{c|}{$3$} & 5.16654423E-02 & 4.800859188E-01 & 275\tabularnewline
$4$ & 1.93312527E-02 & 1.275050913E-01 & 650\tabularnewline
$5$ & 6.28902440E-03 & 7.890462696E-02 & 1106\tabularnewline
$6$ & 7.75971037E-04 & 1.230099037E-02 & 2324\tabularnewline
$7$ & 2.26932864E-04 & 1.648187450E-03 & 4211\tabularnewline
$8$ & 3.16373881E-05 & 1.992534346E-04 & 5781\tabularnewline
\hline 
\end{tabular}

\caption{$\|E_{\mathbf{u}}\|_{1}$ and $\|E_{p}\|_{0}$ for various values
of $W$}

\end{table}

So far we have considered the Stokes problem with Dirichlet condition
on the boundary. Here we consider Stokes problem with mixed boundary
conditions. We consider the following Neumann type boundary conditions
on some part of the boundary of the domain \cite{bebeskuc,manou}
\begin{eqnarray*}
\gamma_{_{N}}(\mathbf{u},p)=\frac{\partial\mathbf{u}}{\partial\mathbf{n}}-p\mathbf{n}=\mathbf{g}^{N}\,\,\textrm{or}\\
\gamma_{_{N}}(\mathbf{u},p)=\left(\left(\nabla\mathbf{u}+\nabla\mathbf{u}^{T}\right)-pI\right)\mathbf{n}=\mathbf{g}^{N}
\end{eqnarray*}
 where $\mathbf{n}=(n_{1},n_{2})$ is unit outward normal vector and
$I$ is $2\times2$ identity matrix. 

Details of the existence and regularity of the solution of the Stokes
problem with mixed boundary conditions can be found in \cite{bebeskuc}.
The numerical method proposed in this article also works for mixed
boundary conditions. In this case, we add the following term to the
minimizing functional $\mathcal{R}^{L,W}(\mathbf{u},p)$ defined in
(8)
\begin{eqnarray*}
\sum_{\gamma_{s}\subseteq\Gamma^{N}\cap\Omega_{l}}\left\Vert \gamma_{_{N}}(\mathbf{u},p)-\mathbf{g}^{N}\right\Vert _{\frac{1}{2},\gamma_{s}}^{2}
\end{eqnarray*}
 where $\Gamma^{N}$ is the part of the boundary of the domain on
which the Neumann type of boundary condition is specified. 

\subsubsection*{Example 5: Stokes problem with mixed boundary conditions on a square
domain}

Consider the Stokes problem (equation (1-3) with $\alpha=0$ and $\mu=1$)
on $[0,1]^{2}$ with mixed boundary conditions. Dirichlet boundary
condition is considered on the sides $x=0,x=1$ and $y=1$ and the
following Neumann type boundary condition is taken on the side $y=0$
\begin{eqnarray*}
\frac{\partial\mathbf{u}}{\partial\mathbf{n}}-p\mathbf{n}=\mathbf{g}^{N}.
\end{eqnarray*}

Chosen the data such that the exact solution of the problem is 
\begin{eqnarray*}
u_{1}=sin\,\pi x_{1}\,\,sin\,\pi x_{2},\\
u_{2}=cos\,\pi x_{1}\,\,cos\,\pi x_{2},\\
p=150\,(x_{1}-\frac{1}{2})(x_{2}-\frac{1}{2}).
\end{eqnarray*}

The domain $[0,1]^{2}$ is divided into $4$ elements with equal step
size $h=\frac{1}{2}$ in each direction. The relative errors $\|E_{\mathbf{u}}\|_{1}$,
$\|E_{p}\|_{0}$ and $\|E_{c}\|_{0}$ for various values of $W$ are
shown in Table 6. 

\begin{table}[H]
~~~~~~~~~~~~~~~%
\begin{tabular}{ccccc}
\hline 
$W$ & $\|E_{\mathbf{u}}\|_{1}$ & $\|E_{p}\|_{0}$ & $\|E_{c}\|_{0}$ & Iter\tabularnewline
\hline 
2 & 4.113397300E-01 & 1.58492332E-01 & 1.16867122E-00 & 18\tabularnewline
3 & 9.728812767E-02 & 4.46496670E-02 & 1.89027751E-01 & 94\tabularnewline
4 & 1.580428850E-02 & 5.59220475E-03 & 4.28144621E-02 & 177\tabularnewline
5 & 1.518175283E-03 & 4.89795131E-04 & 4.296603984-03 & 328\tabularnewline
6 & 1.585801041E-04 & 4.96285704E-05 & 4.48482140E-04 & 572\tabularnewline
7 & 1.706519975E-05 & 5.07498200E-06 & 5.12235349E-05 & 968\tabularnewline
8 & 1.489101173E-06 & 4.36770285E-07 & 4.43246894E-06 & 1570\tabularnewline
\hline 
\end{tabular}

\caption{Errors $\|E_{\mathbf{u}}\|_{1},$ $\|E_{p}\|_{0}$ and $\|E_{c}\|_{0}$
for different $W$}
\end{table}

\subsubsection*{Example 6: Stokes problem with mixed boundary conditions on an annular
domain}

Consider the Stokes problem on the annular domain which was considered
in the example 3. The following Neumann type boundary condition is
considered on the side $y=0$
\begin{eqnarray*}
\left((\nabla\mathbf{u}+\nabla\mathbf{u}^{T})-pI\right)\mathbf{n}=\mathbf{g}^{N}.
\end{eqnarray*}
Dirichlet boundary conditions are considered on the other parts of
boundary of the annular domain. 

We have considered the same data as in example 3 and the domain is
divided into 4 elements (see fig. 5). Table 7 shows the relative errors
$\|E_{\mathbf{u}}\|_{1}$ and $\|E_{p}\|_{0}$ for various values
of $W$. The error decays quickly and this shows the exponential accuracy
of the numerical method. The iteration count is also less compared
to the number of iterations in the example 3. 

\begin{table}[H]
~~~~~~~~~~~~~~~~~~~~~~%
\begin{tabular}{cccc}
\hline 
$W$ & $\|E_{\mathbf{u}}\|_{1}$ & $\|E_{p}\|_{0}$ & Iter\tabularnewline
\hline 
2 & 1.26290947E-01 & 8.027834607E-01 & 9\tabularnewline
3 & 6.83953499E-02 & 2.414437390E-01 & 64\tabularnewline
4 & 7.82886019E-03 & 3.957910344E-02 & 124\tabularnewline
5 & 2.84134327E-03 & 7.704636202E-03 & 252\tabularnewline
6 & 9.23139284E-05 & 4.655347741E-04 & 452\tabularnewline
7 & 4.62497985E-05 & 1.716386607E-04 & 851\tabularnewline
8 & 3.99904145E-06 & 1.361305052E-05 & 1063\tabularnewline
9 & 6.49902194E-07 & 2.398434401E-06 & 1642\tabularnewline
10 & 1.26666857E-07 & 3.200850628E-07 & 2314\tabularnewline
\hline 
\end{tabular}

\caption{$\|E_{\mathbf{u}}\|_{1}$ and $\|E_{p}\|_{0}$ for various values
of $W$}
\end{table}

\subsubsection*{Example 7. Generalized Stokes equations on a cube}

Consider the generalized Stokes problem (problem (1-3) in $\mathbb{R}^{3}$
with $\alpha=1$) on the domain ${[-1,1]}^{3}$ with Dirichlet boundary
condition on the boundary. Let $x=(x_{1},x_{2},x_{3})$ be a point
in the domain and $\mathbf{u}=(u_{1},u_{2},u_{3})$ be the velocity
vector. The force function and boundary data are chosen such that
the exact solution of the given problem is given by 
\begin{align*}
u_{1}(x_{1},x_{2},x_{3}) & =4x_{1}^{2}x_{2}x_{3}(1-x_{1})^{2}(1-x_{2})(1-x_{3})(x_{3}-x_{2})\\
u_{2}(x_{1},x_{2},x_{3}) & =4x_{1}x_{2}^{2}x_{2}x_{3}(1-x_{1})(1-x_{2})^{2}(1-x_{3})(x_{1}-x_{3})\\
u_{3}(x_{1},x_{2},x_{3}) & =4x_{1}x_{2}x_{3}^{2}(1-x_{1})(1-x_{2})(1-x_{3})^{2}(x_{2}-x_{1})\\
p(x_{1},x_{2},x_{3}) & =-2x_{1}x_{2}x_{3}+x_{1}^{2}+x_{2}^{2}+x_{3}^{2}+x_{1}x_{2}+x_{1}x_{3}+x_{2}x_{3}-x_{1}-x_{2}-x_{3}.
\end{align*}

Only one element is considered (i.e $[-1,1]^{3}$) and obtained the
approximate solution of the Generalized Stokes problem for $\nu=1,10.$
The table 8 shows the errors $\|E_{\mathbf{u}}\|_{1},$ $\|E_{p}\|_{0}$
and $\|E_{c}\|_{0}$ against different values of $W$ for $\nu=1$
and table 9 shows for $\nu=10.$ 

\begin{table}[H]
\label{tab5} 
\begin{centering}
\begin{tabular}{lcccc}
\hline 
W  & $||E_{u}||_{1}$  & $||E_{p}||_{0}$  & $||E_{c}||_{0}$  & Iter\tabularnewline
\hline 
2  & 1.1458E+01  & 1.5106E+02  & 7.4993E-01  & 41\tabularnewline
4  & 4.3799E-03  & 8.9338E-02  & 1.8608E-02  & 142\tabularnewline
6  & 2.6778E-04  & 4.5485E-03  & 1.4657E-04  & 331\tabularnewline
8  & 3.3323E-06  & 7.1232E-07  & 1.7715E-07  & 1068\tabularnewline
10  & 3.4672E-08  & 6.5445E-08  & 2.4394E-08  & 1776\tabularnewline
\hline 
\end{tabular}
\par\end{centering}
\caption{Errors $\|E_{\mathbf{u}}\|_{1},$ $\|E_{p}\|_{0}$ and $\|E_{c}\|_{0}$
against $W$ for $\nu=1$}
\end{table}

\begin{table}[H]
\label{tab5} 
\begin{centering}
\begin{tabular}{lcccc}
\hline 
W  & $||E_{u}||_{1}$  & $||E_{p}||_{0}$  & $||E_{c}||_{0}$  & Iter\tabularnewline
\hline 
2  & 1.4913E+01  & 1.4505E+01  & 2.9320E-01  & 68\tabularnewline
4  & 3.8626E-04  & 9.5183E-04  & 2.7758E-02  & 279\tabularnewline
6  & 1.1037E-04  & 6.9155E-04  & 8.0836E-05  & 800\tabularnewline
8  & 3.1704E-06  & 7.7042E-07  & 1.8919E-07  & 1797\tabularnewline
10  & 3.3793E-08  & 6.4958E-08  & 2.2812E-08  & 2861\tabularnewline
\hline 
\end{tabular}
\par\end{centering}
\caption{Errors $\|E_{\mathbf{u}}\|_{1},$ $\|E_{p}\|_{0}$ and $\|E_{c}\|_{0}$
against $W$ for $\nu=10$}
\end{table}

The results shows that the errors $\|E_{\mathbf{u}}\|_{1},$ $\|E_{p}\|_{0}$
and $\|E_{c}\|_{0}$ decays exponentially. We have presented the error
in the continuity equation $\left\Vert E_{c}\right\Vert _{0}$ against
$W$ in few other examples also in this section. The decay of $\left\Vert E_{c}\right\Vert _{0}$
shows the mass conservation property of the method. Similar behavior
has been observed in all the other examples too. 

\section{Conclusion and future work}

In this article we have studied the performance of the nonconforming
least-squares spectral element method for Stokes problems on smooth
domains. The generalized Stokes equation, Stokes problem with mixed
boundary conditions and also the Stokes problem on curvilinear domains
were considered. Spectral approximation is nonconforming and same
order spectral element functions are used for both velocity and pressure
variables. The numerical results shows that the method is exponentially
accurate in both $\mathbf{u}$ and $p$. Since the numerical method
is least-squares, the obtained linear system is symmetric positive
definite. In addition to these advantages, the numerical scheme have
good mass conservation property. The decay of the error in continuity
equation in $L^{2}$ norm shows that the method works very well while
conserving the mass. Studying the performance of this approach for
unsteady flow problems on curvilinear domains is under progress. The
Stokes interface problems is under consideration for the future work. 

\section*{Acknowledgements}

\end{document}